\documentclass[12pt,draftclsnofoot,onecolumn]{IEEEtran}
\usepackage{graphicx}
\usepackage{graphics}
\usepackage[dvips]{epsfig}
\usepackage{latexsym}

\newcommand{\R}{I \! \! R}

\newcommand{\C}{I \! \! \! \! {C}}

\newcommand{\bd}{{\bf d}}
\newcommand{\bud}{{\bf \underline{d}}}

\newcommand{\bux}{{\bf \underline{x}}}

\newcommand{\buy}{{\bf \underline{y}}}
\newcommand{\by}{{\bf y}}

\newcommand{\beps}{{\mbox{\boldmath $\underline{\epsilon}$}}}
\newcommand{\buet}{{\mbox{\boldmath $\underline{\eta}$}}}
\newcommand{\buxi}{{\mbox{\boldmath $\underline{\xi}$}}}

\newcommand{\ux}{{\underline x}}
\newcommand{\uy}{{\underline y}}
\newcommand{\uz}{{\underline z}}

\newcommand{\bub}{{\bf{\underline b}}}

\newcommand{\ue}{{\underline e}}

\newcommand{\uxi}{{\underline{\xi}}}

\newcommand{\ub}{{\underline b}}
\newcommand{\ud}{{\underline d}}
\newcommand{\bA}{{\bf A}}

\newcommand{\bQ}{{\bf Q}}

\newtheorem{teo}{Theorem}
\newtheorem{lem}{Lemma}
\newtheorem{prop}{Proposition}
\newcommand{\bb}{\begin{eqnarray*}}
\newcommand{\be}{\end{eqnarray*}}
\newcommand{\bbn}{\begin{eqnarray}}
\newcommand{\ben}{\end{eqnarray}}

\begin{document}
\title{On a Class of  Parameters Estimators in Linear Models Dominating the Least Squares one, Based on   Compressed Sensing Techniques}
\author{Piero Barone, Isabella Lari
 \thanks{P. Barone is with Istituto per le Applicazioni del Calcolo ''M. Picone'',
C.N.R.,
Via dei Taurini 19, 00185 Rome, Italy,
p.barone@iac.cnr.it, piero.barone@gmail.com}%
\thanks{I. Lari is with Dipartimento di Scienze  Statistiche, Universita' Sapienza,
Piazzale Aldo Moro 5, 00185 Rome, Italy,
isabella.lari@uniroma1.it}}
\maketitle

\begin{abstract}
The estimation of parameters in a  linear model is considered under the hypothesis that the noise, with finite second order statistics, can be represented in a given deterministic basis by random coefficients. An extended underdetermined design matrix is then considered and an estimator of the extended parameters is proposed with minimum $l_1$ norm. It is proved that if the noise variance is larger than a threshold, which depends on the unknown parameters and on the extended design matrix, then the proposed estimator of the original parameters dominates the least-squares estimator in the sense of the mean square error. A small simulation illustrates the behavior of the proposed estimator. Moreover it is shown experimentally that the proposed estimator can be convenient even if the design matrix is not known but only an estimate can be used. Furthermore the noise basis can eventually be used to introduce some prior information in the estimation process. These points are illustrated by simulation by using the proposed estimator  for solving a difficult inverse ill-posed problem related to the complex moments of an atomic complex measure.
\end{abstract}



\section*{Introduction}

Linear models are ubiquitous in applied sciences. Parameters estimation methods have been developed since a long time ago.
In order to motivate the approach that we are proposing, we make
  some  considerations on parameters estimation in linear models related to our purpose. Denoting  random quantities by  bold characters, let us consider the model
\bbn \bud=V\uxi+\beps,\,\;V\in \C^{n\times p},\; n\ge p, \;\mbox{rank}(V)=p,\;\uxi\in \C^p \label{eq1} \ben
where $\bud$ is an $n-$variate complex random vector representing the measured data and $\beps$ is a $n-$variate zero mean complex random vector with finite  second moments representing the measuring error.
The design matrix $V$ is assumed to be ill-conditioned w.r. to the inversion i.e. the ratio of its largest to the smallest singular value is large.

In many applications the  parameters vector $\uxi$ to be estimated represents some well defined object about which much a priori information is available. This motivated the introduction of regularization methods which enforces a parameter estimate with expected properties by solving a modified problem e.g. of the form
\bb \hat{\uxi}=\mbox{argmin}_{\ux}\left( \|\ud-V\ux\|^2+\lambda f(\ux) \right)\be
where $f(\ux)\ge 0$ is a regularitazion function which represents the  prior information and $\lambda>0$ is an hyperparameter balancing the fit to the data and the prior information. In a stochastic environment the Bayes paradigm implements the same idea in a more general form. Given a prior distribution of the parameters and a likelihood, a function proportional to the posterior is used to get estimators either by solving an optimization problem or simply by sampling from the posterior. This last approach is able to cope with problems of huge dimension (MCMC). The main advantage of the regularization and the Bayesian approaches is to strongly reduce the ill-conditioning of the problem i.e. different realizations of the data produce essentially the same parameters estimate.

However in many cases the regularization approach makes no sense because no specific prior information is available on the parameters. In these cases the classical least squares estimator is
\bbn \buxi_{LS}=V^+\bud=(V'V)^{-1}V'\bud \label{eq2} \ben
where prime denotes transposition and plus denotes generalized inversion. The first and second order statistics of $\buxi_{LS}$ are
$$E[\buxi_{LS}]=\uxi,\;\; \mbox{cov}[\buxi_{LS}]=V^+\mbox{cov}[\beps]V^{+'}.$$
Therefore the least squares estimator is not distorted and its covariance and mean square error (MSE) reduces to
$$MSE_{LS}=\mbox{tr}(\mbox{cov}[\buxi_{LS}])=\sigma^2\mbox{tr}({(V'V)^{-1}})$$ when the error is identically distributed with variance $\sigma^2$.
Without loss of generality in the following this setup will be assumed.  When the design matrix is ill-conditioned $MSE_{LS}$ can be quite large. In many instances this can be a serious problem because of the consequent instability of the estimates. It is therefore reasonable to allow some bias in the estimators in order to reduce their variability measured by the MSE. Several methods are reported in the literature (e.g.\cite{eldar}) which modify the least-squares estimator according to some criterium. In this work a different approach is pursued with the same aim. The basic observation is that sometimes we are not able to characterize the parameters but we are able to characterize the noise quite well. As an example of this situation we quote the complex exponentials approximation problem \cite{barja2,bdsp} where it is well known that under a suitable coordinate transformation the noise clusters around the unit circle in the complex plane, but for some gaps, accordingly to an equilibrium measure induced by a logarithmic potential \cite{barja}. Moreover this behavior is quite general (universal) i.e. it does not depend on the specific distribution of the noise \cite{brmta}. The idea is then to consider a model for the noise
\bbn \beps=V_e\buet,\;\;V_e\in \C^{n\times (m-p)},\;m>p \label{eq3}\ben
where the matrix $V_e$ is assigned on the basis of the assumed information about the noise and $\buet$ is a random complex  vector of noise-related parameters to be estimated. We can then consider an extended  model
\bb \bud=V\uxi+V_e\buet=[V | V_e][\uxi'|\buet']'=A\bux,\\ A\in \C^{n\times m},\; n\le m, \;\mbox{rank}(A)\ge p \label{eq4} \be
which is underdetermined if $m>n$. We have now a problem similar to a compressed sensing problem (\cite{don,cand}) with the important simplification that we know which are the noise related components in the extended design matrix $A$.
In order to exploit this similarity we make use of the real isomorph transformation to reformulate the problem in real variables
\bb \ux\in\C^{r}\rightarrow\left[\begin{array}{ll}
\Re{\ux}\\
\Im{\ux}
  \end{array}\right] \in \R^{2r},\\
X  \in\C^{r\times s}\rightarrow\left[\begin{array}{ll}
\Re{X} & -\Im{X}\\
\Im{X} & \;\;\;\Re{X}
  \end{array}\right] \in \R^{2r\times 2s}.
\be
Hence in the following only the real case will be discussed but in the last section where we implicitly assume that the real isomorph transformation is used.
We can then consider the estimator given by
\bbn \left\{\begin{array}{ll}
\hat{\bux}=\mbox{argmin}_{\bux}\|\bux\|_1 \\
\bud=A\bux
  \end{array}\right. \;\;\mbox{ where } \hat{\bux}= \left[\begin{array}{ll}\buxi_D \\
  \hat{\buet}\end{array}\right]\label{eq5}
\ben
and find conditions on $\uxi$, $V_e$ and $\sigma^2$ such that
\bbn MSE_D=E[\|\uxi-\buxi_D\|_2^2] < MSE_{LS}.\label{eq6}\ben
We notice that, by introducing Lagrange multipliers, the problem above is equivalent to a regularization problem with a special regularization function given by the $l_1$ norm of the extended unknown vector.
In section one an explicit form of the estimator is provided. In section two conditions on $\uxi$, $V_e$ and $\sigma^2$ are derived. In section three a small simulation related to the difficult problem of complex exponential approximation is performed to illustrate the advantages of the proposed estimator.

\section{Explicit form of the estimator}

In order to get an explicit form of the estimator $\buxi_D$ let us consider the case $m=n+p$. We first state the following
\begin{lem}
If $B\in \R^{n\times p}$ has rank $p$, the problem
\begin{equation} \label{Problem1}
z^*=\min_{\ux}\sum_{i=1}^n|(\ub-B\ux)_i|
\end{equation}
has at least one solution of the form
$$ \hat{\ux}=B_p^{-1}\ub_p $$
where $B_p$ is a non-singular submatrix of $B$ of order $p$ and $\ub_p$ is the corresponding subvector of $\ub$.
\label{lem1}
\end{lem}
\underline{proof.}

\noindent
Any optimal solution $\ux^*$ to problem (\ref{Problem1}) induces a partition $ M^+, M^- $ of $\{1, \ldots,n\}$ such that:
$$
\begin{array}{ll}
(\ub-B\ux^*)_i \geq 0 & i \in M^+\\
(-\ub+B\ux^*)_i > 0 & i \in M^-.
\end{array}
$$

Consider the following Linear Program based on the partition $ M^+, M^- $:
\bbn \label{Problem2}
\begin{array}{ll}
w^*=& \min_{\ux}\sum_{i \in M^+} (\ub-B\ux)_i + \sum_{i \in M^-} (-\ub+B\ux)_i \\
& (\ub-B\ux)_i \geq 0,\;\;\;\; i \in M^+\\
& (-\ub+B\ux)_i \geq 0,\;\;\;\;i \in M^-.
\end{array}
\ben
Problem (\ref{Problem1}) is a relaxation of problem (\ref{Problem2}) and, in particular, for all feasible solutions of (\ref{Problem2}), the objective function of (\ref{Problem1}) is equal to the objective function of (\ref{Problem2}); furthermore, the optimal solution $\ux^*$ to (\ref{Problem1}) is feasible  for (\ref{Problem2}). It follows that $\ux^*$ is optimal also for (\ref{Problem2}) and $z^* = w^*$. Hence, any optimal solution to (\ref{Problem2}) is optimal also for (\ref{Problem1}).\\
Since problem (\ref{Problem2}) has a finite optimal solution and $B$ has rank $p$ then, for the fundamental theorem of Linear Programming, there exists at least an optimal basic feasible solution to (\ref{Problem2}), i.e. an optimal solution of the form
$$ \hat{\ux}=B_p^{-1}\ub_p $$
where $B_p$ is a non-singular submatrix of $B$ of order $p$ and $\ub_p$ is the corresponding subvector of $\ub$.
Hence, by the above results, also problem (\ref{Problem1}) has at least an optimal solution having this form.
$\;\;\;\;\Box$

The following proposition  holds
\begin{prop}
Let $V_e\in \R^{n\times n}$ be a non-singular matrix, then the estimator $\buxi_D$ is given by
\bbn\buxi_D=B_p^{-1}\bub_p \label{eq7} \ben
where $B_p$ is a non-singular  submatrix of order $p$ of the matrix
$$ B=  \left[\begin{array}{l}
-I_{p} \\
V_e^{-1}V
  \end{array}\right]\in\R^{(n+p)\times p}$$
and $\bub_p$ is the corresponding subvector of the vector $\bub= \left[\begin{array}{l}
\underline{0} \\
V_e^{-1}\bud
  \end{array}\right]\in\R^{(n+p)}$.

\end{prop}
\underline{proof.}

\noindent As the matrix $V_e$ is square and non-singular, we can  solve for  $\buet$ equation (\ref{eq4}) getting
\bbn\buet=\buet(\uxi)=V_e^{-1}(\bud-V\uxi).\label{eq8}\ben
But then
\bb \bux &=& \left[\begin{array}{l}
\uxi \\
\buet
  \end{array}\right]=\left[\begin{array}{l}
\uxi \\
V_e^{-1}\bud-V_e^{-1}V\uxi
  \end{array}\right]\\& =&
 \left[\begin{array}{l}
\underline{0} \\
V_e^{-1}\bud
  \end{array}\right]
  -
   \left[\begin{array}{l}
-I_p \\
V_e^{-1}V
  \end{array}\right]\uxi
  \label{eq9}
\be
or
\bbn \bux=\bub-B\uxi,\;\;B\in\R^{(n+p)\times p},\;\;\mbox{rank}(B)=p.\ben
Therefore
$ \|\bux\|_1=\sum_{i=1}^m|(\bub-B\uxi)_i|$ and eq. (\ref{eq5}) becomes
\bbn
\left\{\begin{array}{ll}\buxi_D=\mbox{argmin}_{\uy}\sum_{i=1}^m|(\bub-B\uy)_i| \\
\hat{\buet}=\buet(\buxi_D)\end{array}\right. .\label{eq10}
\ben
By Lemma \ref{lem1} there exists at least one solution of the form
$$ \buxi_D=B_p^{-1}\bub_p $$
where $B_p$ is a non-singular  submatrix of $B$ of order $p$ and $\bub_p$ is the corresponding subvector of $\bub.\;\;\;\Box$

\section{Conditions on $\uxi$, $V_e$ and $\sigma^2$}

We start by studying the simple case where $V$ is made up by the first $p$ columns of the identity matrix $I_n$ and $V_e=I_n.$
The following proposition holds
\begin{prop}
If
$V=[\ue_1,\dots,\ue_p]\in\R^{n\times p}$ and $V_e=I_n$
the estimator given in eq. (\ref{eq7}) dominates the least square estimator if
$$\sigma^2>\frac{\uxi'\uxi}{p}.$$
\end{prop}
\underline{proof.}

\noindent
By hypothesis
$$ \bub=\left[\begin{array}{l}
\underline{0}_p \\
\bud
  \end{array}\right]\in \R^{n+p},
 \;\;
 B=  \left[\begin{array}{ll}
-I_p \\
\;\;I_p\\
0_{(n-p)\times p}
  \end{array}\right]\in\R^{(n+p)\times p}$$
but then
$$\sum_{i=1}^{n+p}|(\bub-B \buy)_i|=\sum_{i=1}^p| \by_i|+\sum_{i=1}^p|\bd_i-\by_i|+\sum_{i=p+1}^n|\bd_i|$$
If $Q=[\ue_{i_1},\dots,\ue_{i_{p}}]'\in \R^{p\times (n+p)}$ is such that $B_p=QB$ is non-singular then
$|QB|$ must be a permutation matrix $P\in \R^{p\times p}$. But then
$$ \by_j=\ue_j'B_p^{-1}\bub_p=\left\{ \begin{array}{l}
0 \\
\pm \bd_{i_j},\;i_j\in\{1,\dots,p\}
  \end{array}\right.$$
  and therefore the minimum value of $\sum_{i=1}^{n+p}|(\bub-B \buy)_i|$ is equal to
$\sum_{i=1}^n|\bd_i|$ and it is obtained when
\bbn \by_j=(\buxi_D)_j=\ue_j'B_p^{-1}\bub_p=\left\{ \begin{array}{l}
0 \\
\bd_{j}
  \end{array}\right.\label{vinc}\ben
  But this can happen in $u=1+\sum_{k=1}^p\left(\begin{array}{l}
p \\
k
  \end{array}\right)$ different ways.
  Denoting by $\mathcal{ I}$ the set of indices $\{i_1,\dots,i_p\}$  which satisfy the constraint given in (\ref{vinc}) we have that $\buxi_D=(\bQ B)^{-1}\bQ\bub$ and $\bQ$ has a uniform distribution in the set $\mathcal{ I}$ of cardinality $u$  independently of the distribution of $\bud$.

  \noindent We then have
\bb MSE_D &=& E[\|\uxi-\buxi_D\|_2^2]\\ &=&E[(\uxi-(\bQ B)^{-1}\bQ\bub)'(\uxi-(\bQ B)^{-1}\bQ\bub)]\\
&=& \uxi'\uxi + E[((\bQ B)^{-1}\bQ W\bud)' (\bQ B)^{-1}\bQ W\bud]-\\ &&2E[((\bQ B)^{-1}\bQ W\bud)']\uxi \be
where $W=\left[\begin{array}{l}
\underline{0}_{p\times n} \\
I_n
  \end{array}\right]$.
Let be $\bA=(\bQ B)^{-1}\bQ W$.
 We notice that $\bA=[\bA_1 | 0_{p\times (n+p)}]$ and $\bA_1$ is a $p\times p$ matrix which is zero everywhere but in the main diagonal where there is a one in the $j-$th row iff $(\buxi_D)_j\ne 0$. Therefore $\bA_1$ is symmetric and idempotent and we have
 \bb MSE_D&=&\uxi'\uxi + E[\bud_p'\bA_1 \bud_p]-2E[\bud_p'\bA_1]\uxi \\&=&\uxi'\uxi +\mbox{tr}\{E[\bud_p\bud_p'\bA_1 ]\}-2E[\bud_p'\bA_1]\uxi\be
 where $\bud_p$ is the restriction of $\bud$ to its  first $p$ components.
 As the distribution of $\bA_1$ is uniform on a finite set with probability $\frac{1}{u}$ of each event independently of the distribution of $\bud$,  we have
 $$E[\bud_p'\bA_1]\uxi=E[\bud_p']\cdot E[\bA_1]\uxi=\frac{\uxi'\uxi}{2}$$
 $$E[\bud_p\bud_p'\bA_1 ]=E[\bud_p\bud_p']\cdot E[\bA_1]=\frac{\sigma^2 I_p+\uxi\uxi'}{2} $$
  and then
$$ MSE_D = \frac{\uxi'\uxi+p\sigma^2}{2}.$$
By imposing the condition $MSE_D<MSE_{LS}$ and noticing that in this case $MSE_{LS}=p\sigma^2$ we get the thesis. $\;\;\;\Box$

Let us now consider the case when $V\in \R^{n\times p}$ is generic. The following proposition holds
\begin{prop}
If $\mbox{rank}(V)= p$, the columns of $V_e$ are the left singular vectors of $V$ and the number of  singular values of $V$ greater than one are $q<p$, then
the estimator given in eq. (\ref{eq7}) dominates the least square estimator if
$$\sigma^2>\frac{\sum_{j=q+1}^p \tilde{\xi}_j^2}{\sum_{j=q+1}^p c_j^{-2}},\;\;\;\;\tilde{\uxi}=U_2\uxi$$
where $c_1\ge c_2\ge \dots\ge c_p \ge 0$ are the singular values of $V$ and the columns of $U_2$ are the right singular vectors of $V$.
\end{prop}
\underline{proof.}

\noindent
Let  $V=U_1 D U_2$ be the singular value decomposition of $V$ where $U_1\in\R^{n\times n}$ and $U_2\in\R^{p\times p}$ are orthogonal and \bb D=\left[\begin{array}{l}
D_p \\
0_{(n-p)\times p}
  \end{array}\right]\in \R^{n\times p},\\ D_p=\mbox{diag}[c_1,\dots,c_p],\;c_1\ge c_2\ge \dots\ge c_q > 1. \be Equation (\ref{eq4}) becomes
\bbn \bud=U_1 D U_2\uxi+V_e\buet \label{svd1} \ben
or, by defining $\tilde{\bud}=U'_1\bud$ and $\tilde{\uxi}=U_2\uxi$, without loss of generality we can consider the model
\bbn \tilde{\bud}= D \tilde{\uxi}+U_1'V_e\buet \label{svdm}. \ben
By hypothesis $V_e=U_1$ therefore the model becomes
$$\tilde{\bud}= D \tilde{\uxi}+\buet$$ and
$$ \bub=\left[\begin{array}{l}
\underline{0}_p \\
\tilde{\bud}
  \end{array}\right]\in \R^{n+p},
 \;\;
 B=  \left[\begin{array}{ll}
-I_p \\
\;\;D_p\\
0_{(n-p)\times p}
  \end{array}\right]\in\R^{(n+p)\times p}.$$
  but then
$$\sum_{i=1}^{n+p}|(\bub-B \buy)_i|=\sum_{i=1}^p| \by_i|+\sum_{i=1}^p|\tilde{\bd}_i-c_i\by_i|+\sum_{i=p+1}^n|\tilde{\bd}_i|$$
and a value of $\buy$ which minimizes this expression is given by
\bbn \by_j=(\tilde{\buxi}_D)_j=\ue_j'B_p^{-1}\bub_p=\left\{ \begin{array}{l}
0 \;\;\;\; \mbox{  if } c_j\le 1\\
\frac{\tilde{\bd}_{j}}{c_j}\;\; \mbox{ if } c_j>1
  \end{array}\right.\label{vinc1}.\ben
If $q\le p$ is the number of $c_j > 1$  and
\bb A&=&\mbox{diag}\left[\frac{1}{c_1},\dots,\frac{1}{c_q},0,\dots,0\right]\\&=&\left[\begin{array}{ll} A_1 &0_{q\times (n-q)}\\ 0_{(p-q)\times q} & 0_{(p-q)\times (n-q)} \end{array}\right]\in \R^{p\times n},\;\;A_1\in\R^{q\times q}\be
then $\tilde{\buxi}_D=\left[\begin{array}{l}A_1\tilde{\bud}_q \\ 0_{p-q} \end{array}\right]$ and
\bb MSE_D&=&\tilde{\uxi}'\tilde{\uxi}_q + E[\tilde{\bud}_q'A_1^2 \tilde{\bud}_q]-2E[\tilde{\bud}_q'A_1]\tilde{\uxi} \\ &=&\tilde{\uxi}'\tilde{\uxi} +\mbox{tr}\{E[\tilde{\bud}_q\tilde{\bud}_q']A_1^2 \}-2E[\tilde{\bud}_q']A_1\tilde{\uxi}_q
\be
 where $\tilde{\bud}_q$ is the restriction of $\tilde{\bud}$ to its  first $q$ components and $\tilde{\uxi}_q$ is the same for $\tilde{\uxi}$.
 But, if $D_q$ is obtained by putting to zero the last $p-q$ diagonal elements of $D_p$ we have $D_q=A_1^{-1}$ and $$E[\tilde{\bud}_q\tilde{\bud}_q']=\sigma^2 I_q+D_q\tilde{\uxi}_q\tilde{\uxi}_q'D_q \mbox{ and }E[\tilde{\bud}_q]=D_q\tilde{\uxi}_q$$ therefore
 \bb MSE_D&=&\tilde{\uxi}'\tilde{\uxi} +\mbox{tr}\{(\sigma^2 I_q+D_q\tilde{\uxi}_q\tilde{\uxi}_q'D_q)A_1^2 \}-2\tilde{\uxi}_q'D_q A_1\tilde{\uxi}_q
 \\
&=&\tilde{\uxi}'\tilde{\uxi} +\sigma^2\sum_{j=1}^q\frac{1}{c_j^2}-\sum_{j=1}^q\tilde{\xi}_j^2\be
 As $MSE_{LS}=\sigma^2\sum_{j=1}^p\frac{1}{c_j^2}$ and remembering that $\tilde{\uxi}=U_2\uxi$ we have that $MSE_D<MSE_{LS}$ when
 $$\sigma^2>\frac{\sum_{j=q+1}^p \tilde{\xi}_j^2}{\sum_{j=q+1}^p c_j^{-2}}$$
 Finally we notice that in the original variables
 $$\buxi_{LS}=V^+\bud=U_2'D^+U_1'\bud=U_2'D^+\tilde{\bud}=U_2'\tilde{\buxi}_{LS}$$
 and it is easy to check that  the proposed estimator in the original variables is
 $$\buxi_D=U_2'\tilde{\buxi}_D.$$  Therefore we have
$$MSE_{D}(\tilde{\buxi})= E[\|\tilde{\uxi}-\tilde{\buxi}_D\|_2^2]=E[\|\uxi-\buxi_D\|_2^2]=MSE_D(\buxi)$$
because $U_2$ is orthogonal.
As the same is true for the $MSE_{LS}$, this concludes the proof.
  $\;\;\;\Box$

The proposition above can be  generalized to cope with a generic matrix $V_e\in\R^{n\times (m-p)},\;\;n\ge p,\;\;m-p\ge n$. Let us consider the generalized singular value decomposition of the pair $(V,V_e)$  which is given by
$$V=X \mathcal{A} U_1,\;\; V_e=X \mathcal{B} U_2,\;\;X\in\R^{n\times n} \mbox{ invertible }$$ $$U_1\in\R^{p\times p} \mbox{ and }U_2\in\R^{(m-p)\times (m-p)} \mbox{ orthogonal }$$
$$\mathcal{A}'\mathcal{A}+\mathcal{B}'\mathcal{B}=I_n$$ where
\bb \mathcal{A}&=&\left[\begin{array}{l}
0_{(n-p)\times p} \\
\mathcal{A}_p
  \end{array}\right]\in \R^{n\times p},\\ \mathcal{A}_p&=&\mbox{diag}[\alpha_1,\dots,\alpha_p],\;0\le\alpha_1\le \alpha_2\le \dots\le \alpha_p\le 1. \be
\bb \mathcal{B}=\left[\begin{array}{lll}
I_{n-p} & 0_{(n-p)\times p} & 0_{(n-p)\times (m-p-n)}\\
0_{p\times(n-p)} & \mathcal{B}_p & 0_{p\times(m- p-n)}
  \end{array}\right]\in \R^{n\times (m-p)},\be  \bb \mathcal{B}_p=\mbox{diag}[\beta_1,\dots,\beta_p],\;1\ge\beta_1\ge \beta_2\ge \dots\ge \beta_p \ge 0. \be
The following theorem holds
\begin{teo}
If $\mbox{rank}(V)= p\;\;$, $m-p\ge n,\,\;$  $\beta_p>0$ and the number $q$ of ordered pairs $(\alpha_j, \beta_j)$ such that $\alpha_j> \beta_j$ is strictly less than $p$, then
the estimator given in eq. (\ref{eq7}) dominates the least square estimator if
\bbn\sigma^2>\frac{\sum_{j=1}^{p-q} \tilde{\xi}_j^2}{\sum_{j=1}^{p-q} \beta_j^2\alpha_j^{-2}},\;\;\;\;\tilde{\uxi}=U_1\uxi.\label {thresh}\ben

\end{teo}
\underline{proof.}

\noindent
Equation (\ref{eq4}) becomes
\bbn \bud=X \mathcal{A} U_1\uxi+X \mathcal{B} U_2\buet \label{svd} \ben
or, by defining $\tilde{\bud}=X^{-1}\bud$, $\tilde{\uxi}=U_1\uxi$ and $\tilde{\buet}=U_2\buet$, without loss of generality we can consider the model
\bbn \tilde{\bud}= \mathcal{A} \tilde{\uxi}+\mathcal{B} \tilde{\buet} \label{svdm1}. \ben
We then have
$$ \bub=\left[\begin{array}{l}
\underline{0}_p \\
\mathcal{B}^{+}\tilde{\bud}
  \end{array}\right]\in \R^{m},
 \;\;
 B=  \left[\begin{array}{ll}
-I_p \\
\;\;\mathcal{B}^{+}\mathcal{A}
  \end{array}\right]\in\R^{m\times p}$$
  where
  $$\mathcal{B}^{+}=\left[\begin{array}{lll}
I_{n-p} & 0_{(n-p)\times p} \\
0_{p\times(n-p)} & \mathcal{B}_p^{-1} \\
 0_{(m-p-n)\times(n-p) } &  0_{(m- p-n)\times p}
  \end{array}\right]\in \R^{(m-p)\times n}.$$
  Therefore
 $$ \bub=\left[\begin{array}{l}
\underline{0}_p \\
\tilde{\bud}_{n-p} \\
\mathcal{B}_p^{-1}\tilde{\bud}_p\\
0_{(m-p-n)}
  \end{array}\right]$$
  where $\tilde{\bud}_{n-p}=[\tilde{\bd}_1,\dots,\tilde{\bd}_{n-p}],\;\;\tilde{\bud}_p=[\tilde{\bd}_{n-p+1},\dots,\tilde{\bd}_{n}]$ and
 $$ B=  \left[\begin{array}{ll}
-I_p \\
0_{(n-p)\times p}\\
\;\;\mathcal{B}_p^{-1}\mathcal{A}_p\\
0_{(m-p-n)\times p}
  \end{array}\right].$$
  But then
 $$\sum_{i=1}^m|(\bub-B \buy)_i|=\sum_{i=1}^p| \by_i|+\sum_{i=1}^{n-p}|\tilde{\bd}_i|+\sum_{i=1}^p\left|\frac{\tilde{\bd}_{n-p+i}}{\beta_i}-\frac{\alpha_i}{\beta_i}\by_i\right|$$
and a value of $\buy$ which minimizes this expression is given by
\bbn \by_j=(\tilde{\buxi}_D)_j=\ue_j'B_p^{-1}\bub_p=\left\{ \begin{array}{l}
0 \;\;\;\; \;\;\;\;\;\;\mbox{  if } \alpha_j\le \beta_j\\
\frac{\tilde{\bd}_{n-p+j}}{\alpha_j}\;\; \mbox{ if } \alpha_j> \beta_j
  \end{array}\right.\label{vinc2}.\ben
 Noticing that $\alpha_j/ \beta_j$ is an increasing sequence, denoting by  $q\le p$  the number of $\alpha_j> \beta_j$  and if
\bb A&=&\mbox{diag}\left[0,\dots,0,\frac{1}{\alpha_{p-q+1}},\dots,\frac{1}{\alpha_p}\right]\\&=&
\left[\begin{array}{ll} 0_{(p-q)\times (n-q)} &0_{(p-q)\times q}\\ 0_{q\times (n-q)} & A_1 \end{array}\right]\in \R^{p\times n},\;\;A_1\in\R^{q\times q}\be
then $\tilde{\buxi}_D=\left[\begin{array}{l} 0_{p-q} \\ A_1\tilde{\bud}_q \end{array}\right]$
where $\tilde{\bud}_q$ is obtained by taking  the last $q$ components of $\tilde{\bud}$.
It turns out that
\bb MSE_D=\sum_{j=1}^{p-q}\tilde{\xi}_j^2+\sigma^2\sum_{j=p-q+1}^p\frac{\beta_j^2}{\alpha_j^2}\be
\bb MSE_{LS}=\sigma^2\sum_{j=1}^p\frac{\beta_j^2}{\alpha_j^2} \be
and  the proof follows by the same arguments used in the proof of Proposition 2. $\Box$

\noindent\underline{Remark 1}
We notice that the squared bias of the proposed estimator is
$$b^2=\sum_{j=1}^{p-q}\tilde{\xi}_j^2$$
the larger $q$ the smaller $b^2$. The variance is instead controlled by the values of $\beta_j,\;j=p-q+1,\dots,p$. As
$$0\le\alpha_1\le \alpha_2\le \dots\le \alpha_p\le 1$$ and $$1\ge\beta_1\ge \beta_2\ge \dots\ge \beta_p > 0$$
 the best choice to have both bias and variance as small as possible is to choose
$$q=p-1,\;\;\beta_1=\alpha_1,\;\;\beta_j=\epsilon,\;j=2,\dots,p,\;\;0<\epsilon<\alpha_1.$$
We then get
$$b^2=\tilde{\xi}_1^2,\;\;\mbox{ var }=\sigma^2\epsilon^2\sum_{j=2}^p\frac{1}{\alpha_j^2}$$
and the constraint (\ref{thresh}) becomes  $\sigma^2>b^2$. However this constraint can be too strong if the noise is not so large. Therefore it can be convenient to decrease the threshold on $\sigma^2$ by choosing $q<p-1$. In fact we notice that in (\ref{thresh}) the denominator is greater than one and than it can compensate for the larger numerator induced by the choice
$q<p-1$ if $V_e$ is chosen appropriately.

\noindent\underline{Remark 2}
We notice that $$MSE_{LS}-MSE_D=\sigma^2\sum_{j=1}^{p-q}\frac{\beta_j^2}{\alpha_j^2}-b^2$$
is a linear function of $\sigma^2$ whose slope can somewhat be controlled by $\beta_1, \dots, \beta_{p-q}$. If we have an upper bound on the $l_2$ norm of the true parameters vector
$$\tau_b\ge\|\uxi\|_2^2=\|\tilde{\uxi}\|_2^2\ge b^2=\sum_{j=1}^{p-q}\tilde{\xi}_j^2$$ we can not increase the MSE by more than $\tau_b$ by using the proposed method instead than the least squares one, i.e.
$$MSE_{LS}\ge MSE_D-\tau_b .$$  Moreover if we can find $V_e$ such that
$$\sum_{j=1}^{p-q}\frac{\beta_j^2}{\alpha_j^2}=\frac{\tau_b}{\sigma^2}$$
then the proposed method is convenient.

\section{Experimental results}

To illustrate the advantages of the proposed estimator three simulation experiments were performed to compare the distribution of $MSE_D$ and $MSE_{LS}$ in a specific complex exponentials problem. Let us consider the complex model
\bb f(t;p,P)= \sum_{j=1}^p
\xi_j z_j^t,\;\;t\in\R^+\\P=\{\xi_j,z_j,\;j=1,\dots,p\}\in\C^{2p}
\be
 and assume that we want
to estimate $p$ and $P$ from the data
$$\bd_k=f(k\Delta)+\beps_k,\;k=0,\dots,n-1,\;\;\;\Delta>0,\;\;\;n\ge2p$$
with the identifiability condition
$|\mbox{arg}(z_j)|\Delta\le\pi,\;|\mbox{arg}(\xi_j)|\le \pi \;\forall j$, where the noise
$\beps_k$ are i.i.d. zero-mean complex Gaussian
variables with  variance $\sigma^2$ i.e. the real and
imaginary parts of $\bd_k$ are independently distributed as
Gaussian variables with variance $\sigma^2/2$ and mean $\Re
e[f(k\Delta)],\Im m[f(k\Delta)]$ respectively.

The problem arise in many different fields (see e.g.\cite{bdsp} for a short list). It is an inverse problem which can be severely ill posed. In \cite{bdsp} a method is proposed to solve it stably which performs better than standard alternatives. The most difficult part of the problem  is, apparently, the estimation of $p$ and $\uz=[z_1,\dots,z_p]$ because of the non linear dependence of these parameters on the data. The method proposed in \cite{bdsp} concentrates in fact on this part of the problem and solves the Vandermonde linear system in the unknowns $\uxi=[\xi_1,\dots,\xi_p]$, given $p$ and $\uz$, by standard least squares method. However this system  can be very ill posed too because of the large condition number of the Vandermonde matrix if the Euclidean distance $|z_j-z_h|$ of one or more pairs $(z_j,z_h)$ is small.

In the first experiment we assume to know the parameters
$p$ and $\uz$ and we concentrate on the estimation of the parameters $\uxi$.
More precisely we consider the model given in eq.(\ref{eq4}) with $$V(k,h)=z_h^{k-1},\;k=1,\dots,n;\;\;h=1,\dots,p$$  $$V_e(k,h)=e^{\frac{2\pi i(k-1)(h-1)}{m-p}},\;k=1,\dots,n;\;\;h=1,\dots,m-p,$$ $\;m-p\ge n,$
where the choice of $V_e$ is justified by the error model suggested in \cite[sect.1]{barja2}.
The matrix $A=[V|V_e]$ is then scaled as follows
$$\tilde{A}=A \mathcal{D},\;\;\mathcal{D}=\mbox{ diag }[\|A\ue_1\|_2^{-1},\dots,\|A\ue_m\|_2^{-1}]$$
where $\ue_k$ is the $k-$th column of the identity matrix of order $m$ in order to give the same weight to each column of $A$. Best results were obtained for $m=2n$ in this specific case. The following problem is then solved instead than the one with equality constraints given in eq. (\ref{eq5})
\bb \left\{\begin{array}{ll}
\buy=\mbox{argmin}_{\bux}\|\bux\|_1 \\
\|\bud-\tilde{A}\bux\|_2\le \tau
  \end{array}\right. \;\;\mbox{ and }\\ \hat{\bux}=\mathcal{D}\buy =\left[\begin{array}{ll}\buxi_D \\
  \hat{\buet}\end{array}\right],\;\;\;0< \tau \ll \sigma
\be
to cope with eventual numerical not positive definiteness of the matrix $\tilde{A}\tilde{A}'$.
A log-barrier method described in \cite[ch.11]{boyd} is used.
The set of $p=5$ true parameters \bb \uz=\{e^{-0.3-i 2\pi
0.35}, e^{-0.1-i 2\pi  0.3},e^{-0.05-i 2\pi
0.28},\\ e^{-0.0001+i 2\pi 0.2},e^{-0.0001+i 2\pi  0.21}  \}\be $$\uxi=\{20,6,3,2,1\}$$ was considered and
two simulations were performed with variance given respectively by $\sigma^2=100$ and $\sigma^2=200$ which corresponds to a $SNR=0.01$ and $SNR=0.02$ if the signal-to-noise ratio is defined as
 $$SNR=  2\min_j\frac{|\xi_j|^2}{\sigma^2} .$$
 In each simulation $R=1000$ independent realizations of $\bud$ were computed. For each of them the relative errors
$$ E_D(r)= \frac{\|\uxi-\uxi^{(r)}_D\|_2}{\|\uxi\|_2}$$
$$ E_{LS}(r)= \frac{\|\uxi-\uxi^{(r)}_{LS}\|_2}{\|\uxi\|_2}$$
were collected. Their empirical distributions are shown in Fig.1.
We notice that the distribution of $ E_D$ is shifted to the left w.r. to the distribution of $ E_{LS}$ and this effect is more evident for the smallest $SNR$ as expected (see the Remark at the end of the previous section).

In the second experiment  we no longer assume to know $p$ and $\uz$ but instead we estimate them by the method described in \cite{bdsp}. Of course this case no longer fits the theory exposed above because there is a critical font of variability in the design matrix $V$ itself which is very difficult to account for. However we experimentally show that the qualitative results are the same i.e. the proposed method produces estimates better than the least squares ones w.r. to the MSE and the advantage is increasing with the noise variance. In this case we have an estimate  $\hat{p}$ of $p$ and an estimate $\hat{\uz}$ of $\uz$. When $\hat{p}<p$ we can not estimate the whole vector $\uxi$ and the corresponding realization is thrown away. Two simulations are performed with $\sigma^2=4$ and $\sigma^2=2$ and the results are shown in Fig.2.

In the third experiment we notice that the matrix $V_e$ can be used to include some form of prior information on the solution. For example, in the case considered,    $V_e$ is a Vandermonde matrix  based on complex numbers equispaced on the unit circle, while $V$ is a Vandermonde matrix based on the numbers $\uz$. We notice that $z_4$ and $z_5$ are close to the unit circle and therefore they are also close to $e^{\frac{2\pi i(h-1)}{m-p}}$ for some $h$. The information conveyed by the corresponding columns of $V_e$ reinforces the information conveyed by the columns of $V$ associated to $z_4$ and $z_5$. This can be useful  when $\hat{p}<p$. In fact if in the estimation process we consider not only  $\buxi_D$ but also
  $\hat{\buet}$, i.e. the whole vector $\hat{\bux}$, sort its components in decreasing order of their absolute value and consider the first $p$ components as estimates of $\uxi$  we get the results shown in Fig.3. A slight improvement w.r. to the results shown in Fig.2 can be noted.

The mean of the  relative errors of the parameters over the $R=1000$ replications obtained in the three experiments are reported in Table 1.

\begin{table}[tbh]
\begin{scriptsize}
\begin{center}
\begin{tabular}{|c|c|c|c|}
\hline\hline
Experiment&$SNR$&$E_{LS}$&$E_D$\\
\hline
1&0.01&0.62&0.55\\
\hline
1&0.02&0.44&0.41\\
\hline\hline
2&0.5&0.27&0.25\\
\hline
2&1.0&0.21&0.17\\
\hline\hline
3&0.5&0.27&0.21\\
\hline
3&1.0&0.21&0.16\\
\hline\hline
\end{tabular}
\end{center}
\caption{Estimated relative errors in the three experiments for two SNRs. } \label{tb1}
\end{scriptsize}
\end{table}

\section{Conclusions}
It is proved that when estimating the parameters of a  linear model with ill conditioned (w.r. to the inversion) design matrix, it can be convenient to look for a suitable basis for the noise and try the proposed estimator in order to improve the average mean square error of the estimates. Despite of the fact that the proposed method is convenient only if the noise variance is larger than a threshold, which depends on the unknown true parameters vector, it is enough to have an upper bound on its $l_2$ norm to decide if the proposed method is convenient.

\newpage

\begin{figure}
\centering{
\includegraphics[height=5cm,width=7.5cm]{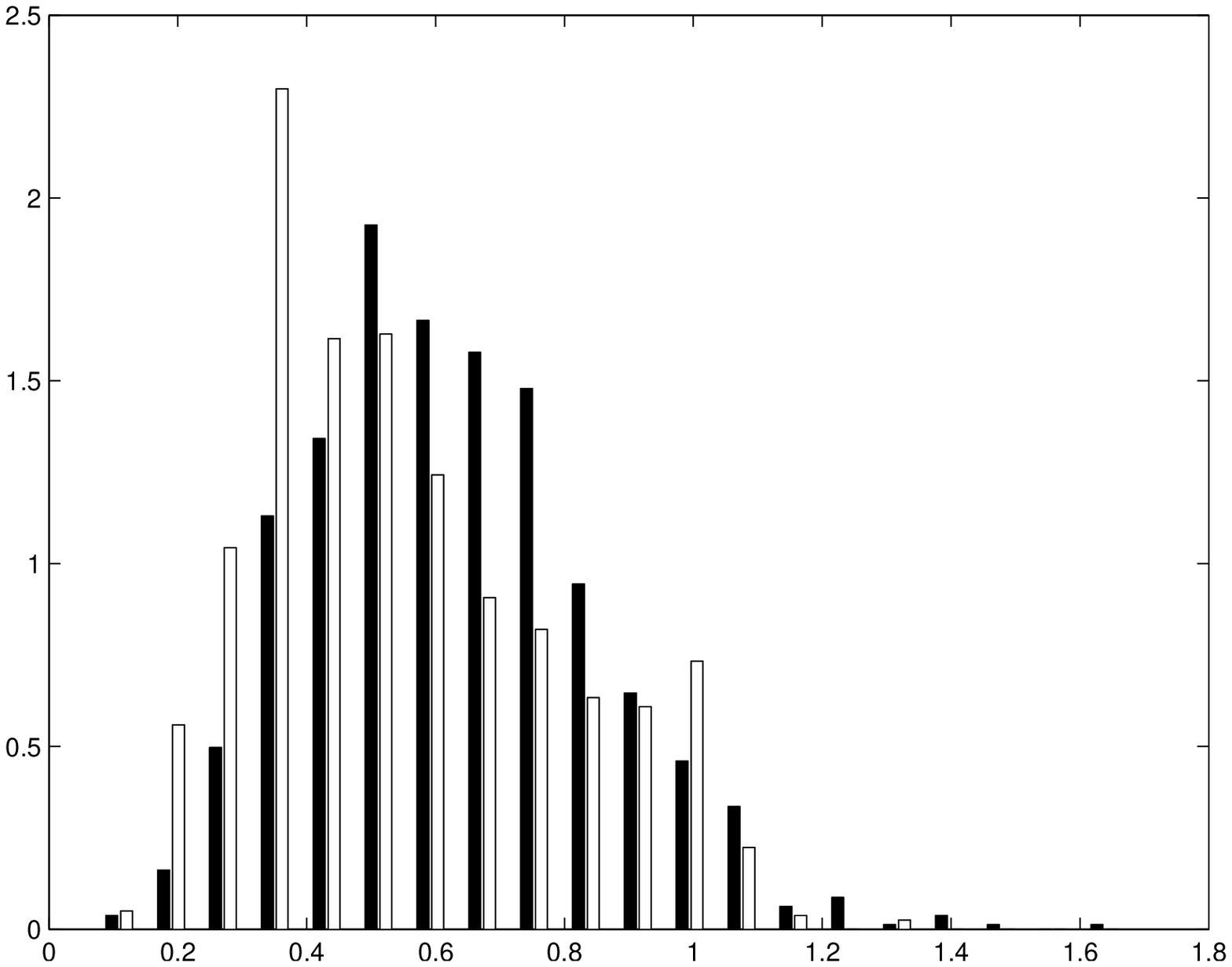}}
\includegraphics[height=5cm,width=7.5cm]{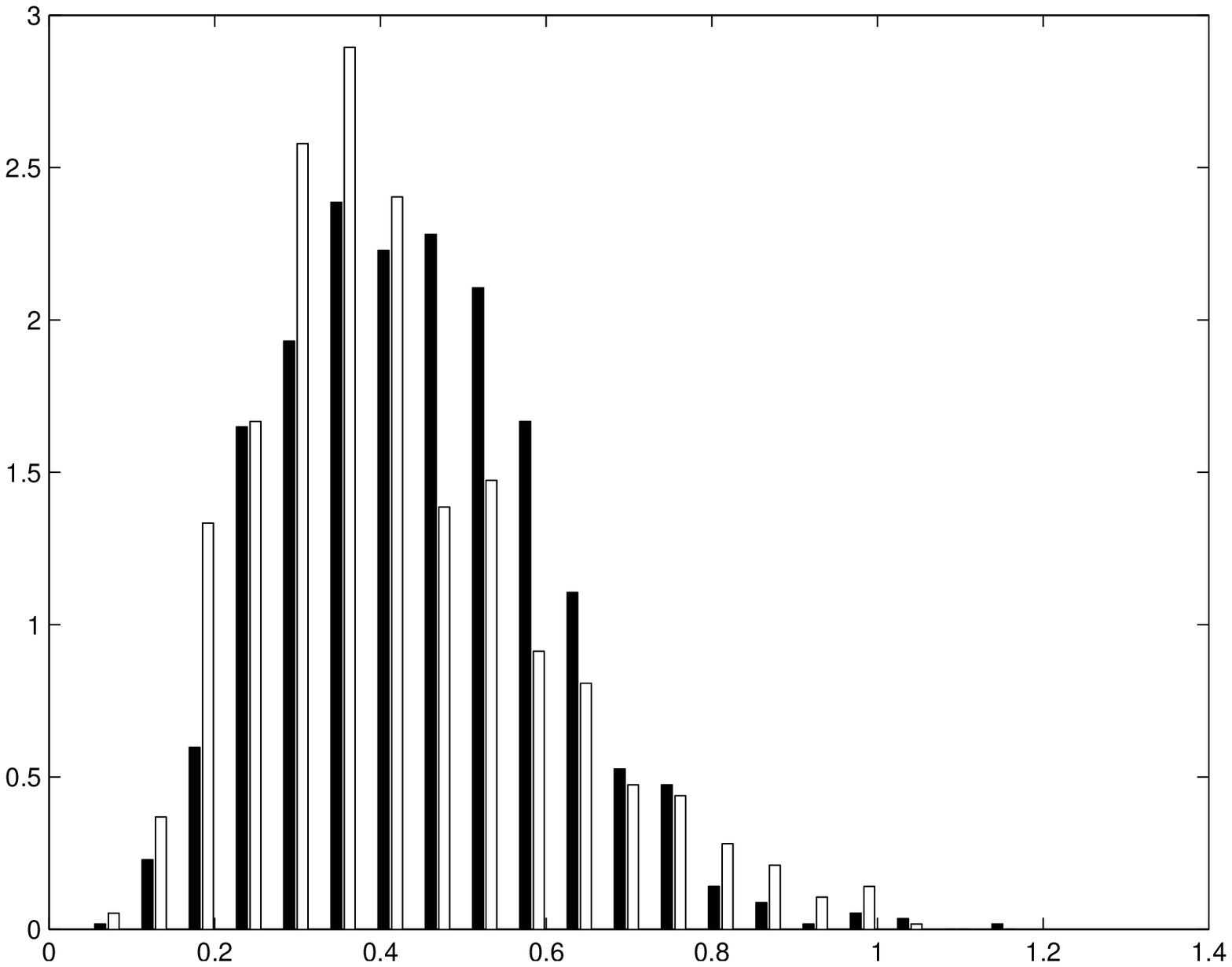}
 \caption{Experiment 1. Distribution of the relative error w.r.to the true parameters in $1000$ replications;  white:  proposed estimates, black: least squares estimates. Left $SNR=0.01$, right $SNR=0.02$. } \label{fig1}
\end{figure}

\begin{figure}
\centering{
\includegraphics[height=5cm,width=7.5cm]{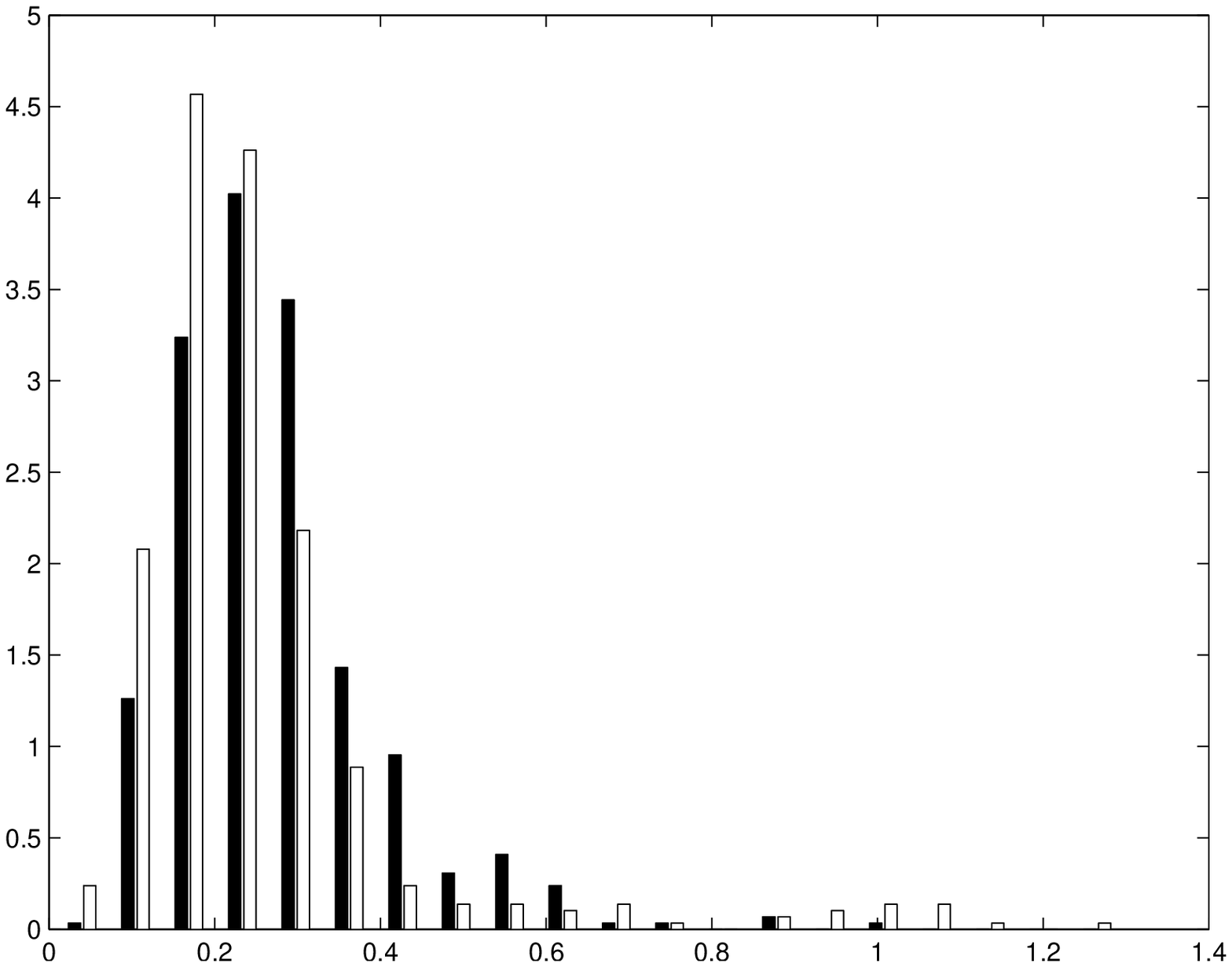}}
\includegraphics[height=5cm,width=7.5cm]{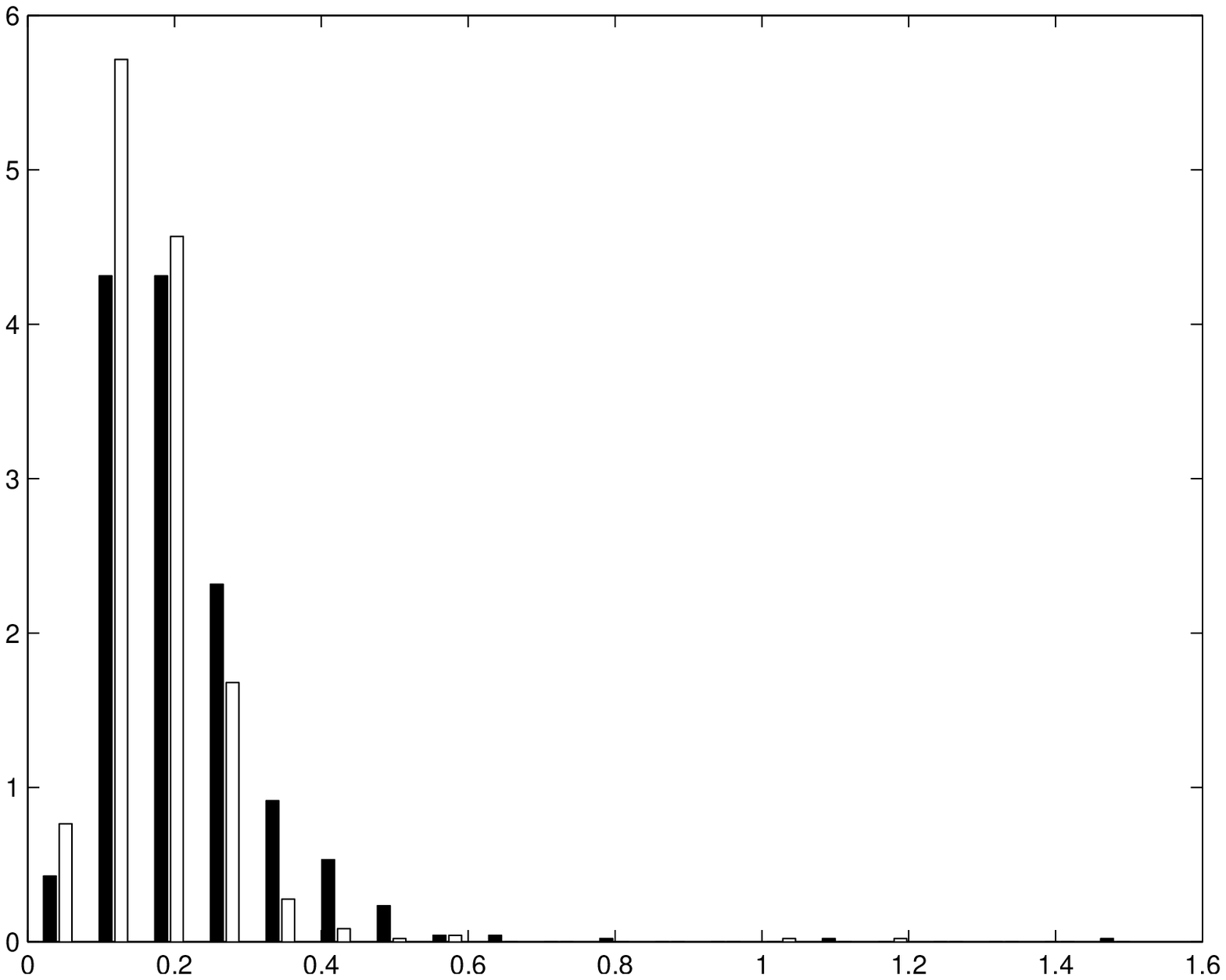}
 \caption{Experiment 2. Distribution of the relative error w.r.to the true parameters in $1000$ replications;  white:  proposed estimates, black: least squares estimates. Left $SNR=0.5$, right $SNR=1$. } \label{fig2}
\end{figure}

\begin{figure}
\centering{
\includegraphics[height=5cm,width=7.5cm]{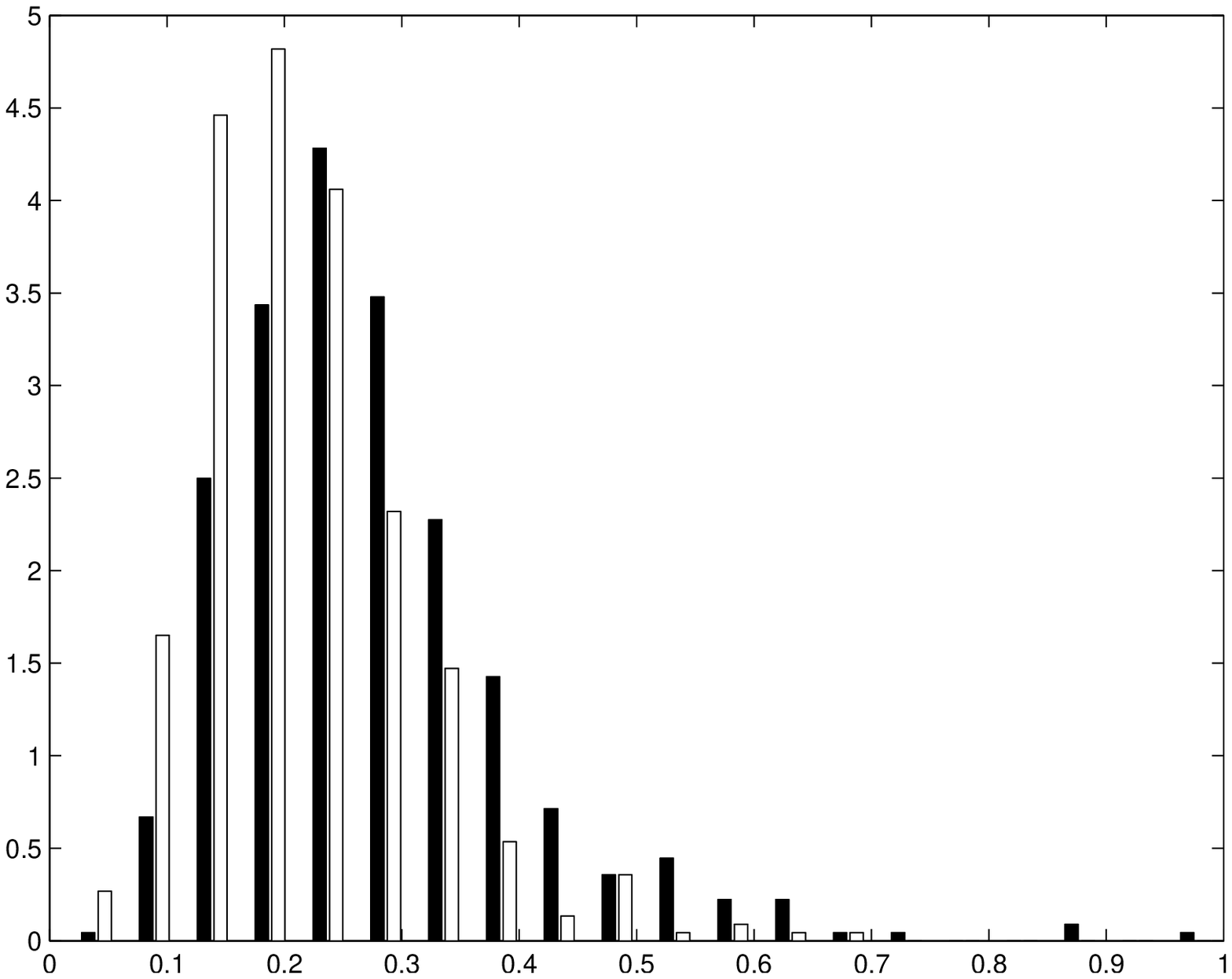}}
\includegraphics[height=5cm,width=7.5cm]{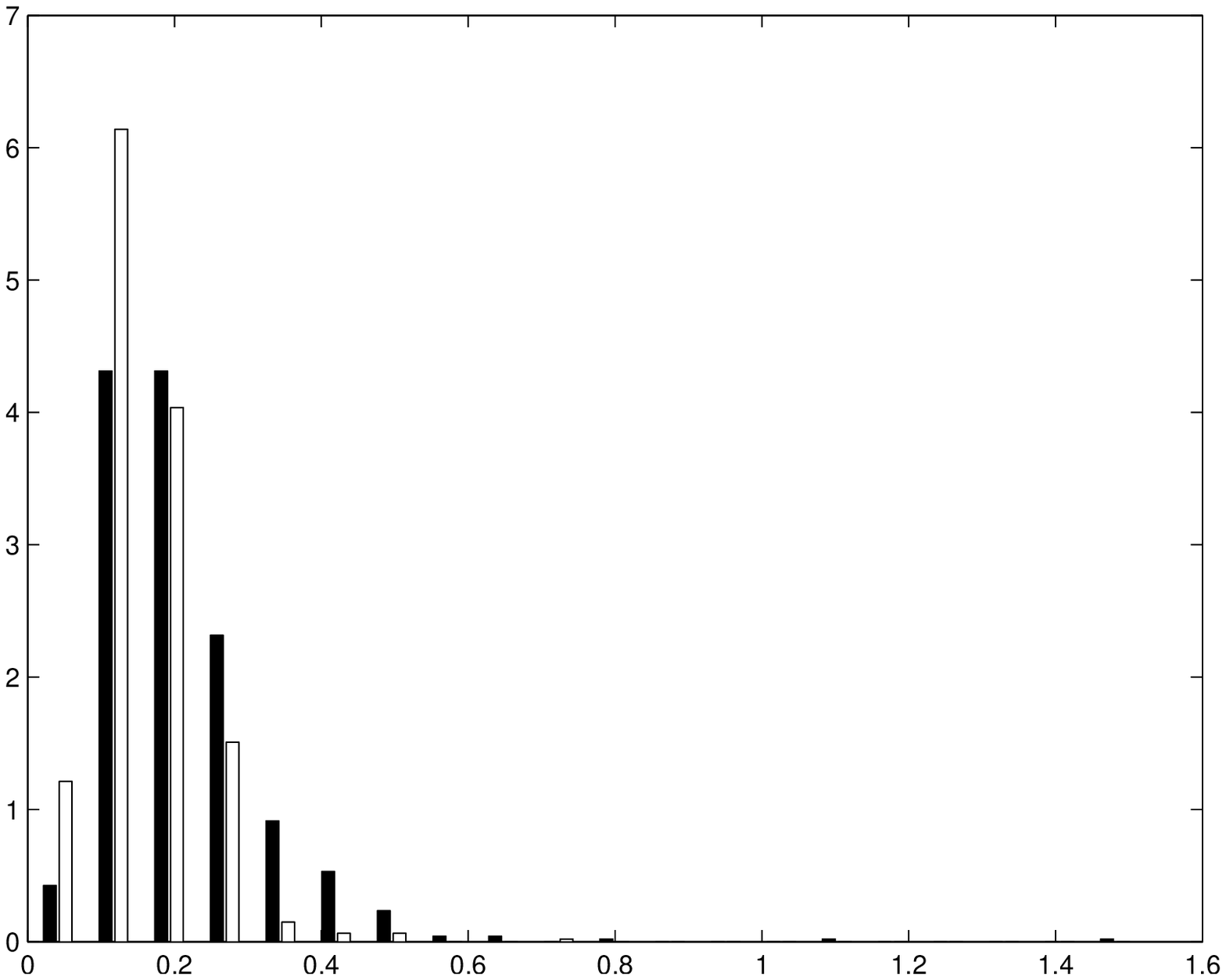}
 \caption{Experiment 3. Distribution of the relative error w.r.to the true parameters in $1000$ replications;  white:  proposed estimates, black: least squares estimates. Left $SNR=0.5$, right $SNR=1$. } \label{fig3}
\end{figure}


\begin{thebibliography}{1}


\bibitem{barja} P. Barone, On the distribution
    of poles of Pade' approximants to the Z-transform
    of complex Gaussian white noise, {\it J. Approx.
    Theory} {\bf 132} (2005) 224--240.
\bibitem{barja2} P. Barone,   A new transform for solving
    the noisy complex exponentials approximation problem, {\it J. Approx.
    Theory} {\bf 155} (2008), 1--27.
\bibitem{brmta} P. Barone, On the universality of the distribution of the generalized eigenvalues of a pencil of Hankel random matrices, {\it Random Matrices: Theory and Applications}, {\bf 2,1} (2013), 1--14, DOI: 10.1142/S2010326312500141
\bibitem{bdsp} P. Barone, A black box method for solving the complex exponentials approximation problem, {\it Digital Signal Processing} {\bf 23} (2013) 49--64
\bibitem{boyd} S. Boyd, L. Vandenberghe, Convex Optimization. Cambridge University Press, 2004.
\bibitem{eldar} Z. Ben-Haim, Y.C. Eldar, Blind minimax estimation, {\it IEEE Trans.Inf.Theory}, {\bf 53,9} (2007), 3145--3157
\bibitem{don} D.L. Donoho, Compressed sensing, {\it IEEE Trans.Inf.Theory}, {\bf 52,4} (2006), 1289--1306
\bibitem{cand} E.J. Candes, J. Romberg, T. Tao, Robust uncertainity principles: exact signal reconstruction from higly incomplete frequency information, {\it IEEE Trans.Inf.Theory}, {\bf 52,2} (2006), 489--509

\end{thebibliography}
\end{document}